\documentclass[12pt,twoside]{paper}
\usepackage{amsmath}
\usepackage{hyperref}
\usepackage{amsfonts}
\usepackage{amssymb}
\usepackage{amsthm}

\newtheorem{theorem}{Theorem}[]

\usepackage{graphics}
\usepackage[pdftex]{graphicx}

\def\be{\begin{displaymath}}
\def\ee{\end{displaymath}}
\def\bee{\begin{equation}}
\def\eee{\end{equation}}
\def\Mas{Ma$\acute{\rm s}$lanka  }
\usepackage{graphicx}

\def\be{\begin{displaymath}}
\def\ee{\end{displaymath}}
\def\bee{\begin{equation}}
\def\eee{\end{equation}}
\def\Mas{Ma$\acute{\rm s}$lanka  }

\begin{document}

\thispagestyle{empty}
\centerline{}
\bigskip
\bigskip
\bigskip
\bigskip
\bigskip
\centerline{\Large\bf Equivalence of Riesz and Baez-Duarte}
\bigskip
\centerline{\Large\bf   criterion for the Riemann Hypothesis}
\bigskip

\begin{center}
{\large \sl J.Cis{\l}o, M. Wolf}\\*[5mm]

Institute of Theoretical Physics, University of Wroc{\l}aw\\
Pl.Maxa Borna 9, PL-50-204 Wroc{\l}aw, Poland\\
\href{mailto:cislo@ift.uni.wroc.pl}{e-mail:cislo@ift.uni.wroc.pl}\\
\href{mailto:mwolf@ift.uni.wroc.pl}{e-mail:mwolf@ift.uni.wroc.pl}\\
\bigskip
\today

\end{center}

\bigskip\bigskip

\begin{center}
{\bf Abstract}\\
\end{center}

\begin{minipage}{12.8cm}
We investigate the relation between the Riesz and  the Baez-Duarte criterion
for the Riemann Hypothesis. In particular we present the relation between the
function $R(x)$ appearing in the Riesz  criterion and the sequence $c_k$ appearing
in the Baez-Duarte formulation. It is shown that $R(x)$ can expressed by $c_k$
and vice versa the sequence $c_k$ can be obtained from the
values of $R(x)$ at integer arguments. We give also some relations involving $c_k$ and $R(x)$,
in particular  value of the alternating sum of $c_k$.
\end{minipage}

\bigskip\bigskip

{\bf 1. Introduction.}\\

The Riemann Hypothesis (RH) states that the nontrivial zeros of the function:
\bee
\zeta(s)=\frac{1}{1-2^{1-s}}\sum_{n=1}^\infty \frac{(-1)^{n-1}}{n^s},
\eee
where $\Re(s)>0$ and $s\neq 1$, are simply and have the real part equal to half, i.e. $\Re (s)=\frac{1}{2}$.
There are probably over 100 statements equivalent to RH, see eg. \cite{Titchmarsh},
\cite{aimath}, \cite{Watkins}.
In the beginning of XX century M. Riesz \cite{Riesz} has considered the function:
\bee
R(x) = \sum_{k=1}^\infty \frac{(-1)^{k+1}x^k}{(k-1)!\zeta(2k)}=
x \sum_{k=0}^\infty \frac{(-1)^{k}x^{k}}{k!\zeta(2k+2)}.
\label{Riesz}
\eee
Unconditionally it can be proved that $R(x)=\mathcal{O}( x^{1/2+\epsilon})$, see
\cite{Titchmarsh}  \S 14.32. Riesz has proved that the Riemann Hypothesis is equivalent
to slower increasing of the function $R(x)$:
\bee
RH \Leftrightarrow   R(x) = \mathcal{O}\left( x^{1/4+\epsilon}\right).
\label{Riesz criterion}
\eee
A few years ago L. Baez-Duarte \cite{Luis2} \cite{Luis3}  considered the sequence
of numbers $c_k$ defined
by:
\bee
c_k=\sum_{j=0}^k {(-1)^j \binom{k}{j}\frac{1}{\zeta(2j+2)}}.
\label{ckmain}
\eee
He proved that RH is equivalent to the following rate of decreasing to zero of the above
sequence:
\bee
RH  \Leftrightarrow c_k={\mathcal{O}}(k^{-\frac{3}{4}+\epsilon})~~~~~~~~~{\rm for~ each~~}  \epsilon>0.
\label{criterion}
\eee
Furthermore, if $\epsilon$ can be put zero, i.e. if
$c_k={\mathcal{O}}(k^{-\frac{3}{4}})$, then the zeros of $\zeta(s)$ are simply. Baez-Duarte
also proved in \cite{Luis3} that it is not possible to replace $\frac{3}{4}$ by larger
exponent. Although  the title of the Baez-Duarte paper was {\it A sequential
Riesz-like criterion for the Riemann Hypothesis} he has not pursued further relation
between $c_k$ and $R(x)$.

In this paper we are going to establish the relation between
$c_k$ and $R(x)$. In  Sect. 2  we will present formulae allowing to obtain values of
$R(x)$ and $c_k$ much faster than from (\ref{Riesz}) and (\ref{ckmain}).
In Sect. 3 we will use the fact that $c_k$ can be obtained as
forward differences of a appropriate sequence to express $R(x)$ in terms of $c_k$. 
Next we will prove equivalence of the Riesz and Baez-Duarte criterion for RH.
In the mathematical logic  the {\it iff} obeys  the transitivity rule:
\be
(p\Leftrightarrow q ~~{\rm AND}~~q \Leftrightarrow s) \Rightarrow (p \Leftrightarrow s)
\ee
thus from (\ref{Riesz criterion}) and (\ref{criterion}) we have that $R(x) = \mathcal{O}\left( x^{1/4+\epsilon}\right)
\Leftrightarrow c_k={\mathcal{O}}(k^{-\frac{3}{4}+\epsilon})$.

However we will prove equivalence (Riesz~ criterion) $ \Leftrightarrow $ (Baez-Duarte~ criterion)
in a more general form, namely the exponents $1/4$ and $3/4$ will be replaced by arbitrary
parameter $\delta$ and  combination $1-\delta$: $c_k=\mathcal{O}\left(k^{-\delta}
\right) \Leftrightarrow R(x)=\mathcal{O}\left(x^{1-\delta}\right)$.
In the final Section we will speculate on some equations involving $c_k$ and $R(x)$,
in particular we will calculate the alternating sum $\sum_{k=0}^\infty (-1)^k c_k$.

\bigskip
\bigskip

{\bf 2. Some facts on the  $R(x)$ and $c_k$}

\bigskip

The most comprehensive source of information about  the Riesz function $R(x)$ we have
found on the Wikipedia
\cite{Wikipedia}. For large negative $x$ function $R(x)$ tends to $xe^{-x}$.
For positive $x$ the behaviour of $R(x)$  is much more difficult to reveal
because the series (\ref{Riesz}) is very slowly convergent. Applying
Kummer's acceleration convergence method  gives
\bee
R(x) = x \sum_{n=1}^\infty \frac{\mu(n)}{n^2} ~e^{-\frac{x}{n^2}}
\label{Riesz2}
\eee

\begin{figure}[pht]
\vspace{-3.5cm}
\begin{minipage}{15.8cm}
\begin{center}
\hspace{-3.5cm}
\includegraphics[width=12cm,angle=0, scale=1]{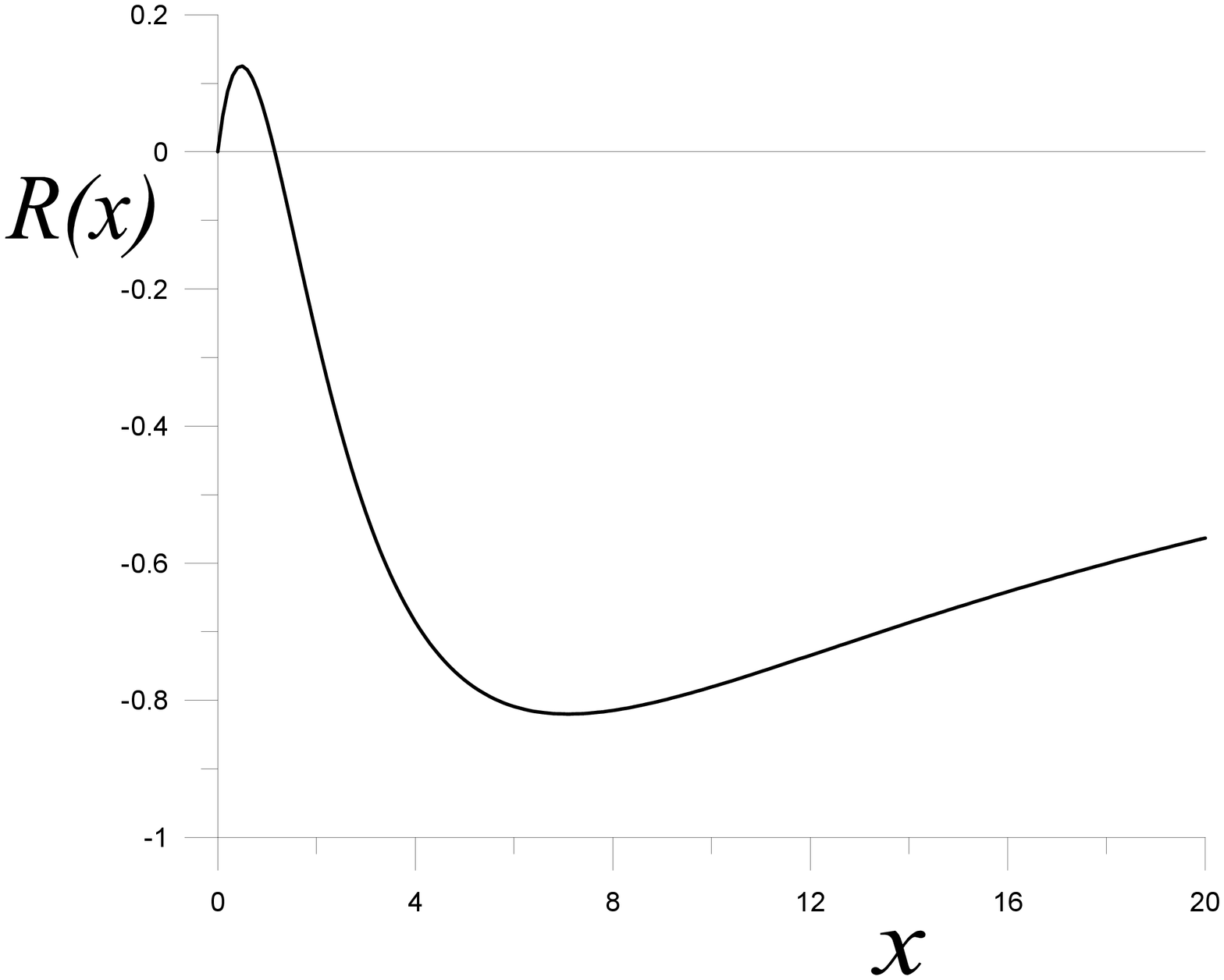} \\
\vspace{-4cm} Fig.1  The plot of $R(x)$ for $x\in(1,20)$. Such a short interval
is chosen to show the first zero of $R(x)$.  \\
\end{center}
\end{minipage}

\vspace{-2.0cm}
\begin{minipage}{15.8cm}
\begin{center}
\hspace{-3.5cm}
\includegraphics[width=12cm,angle=0, scale=1]{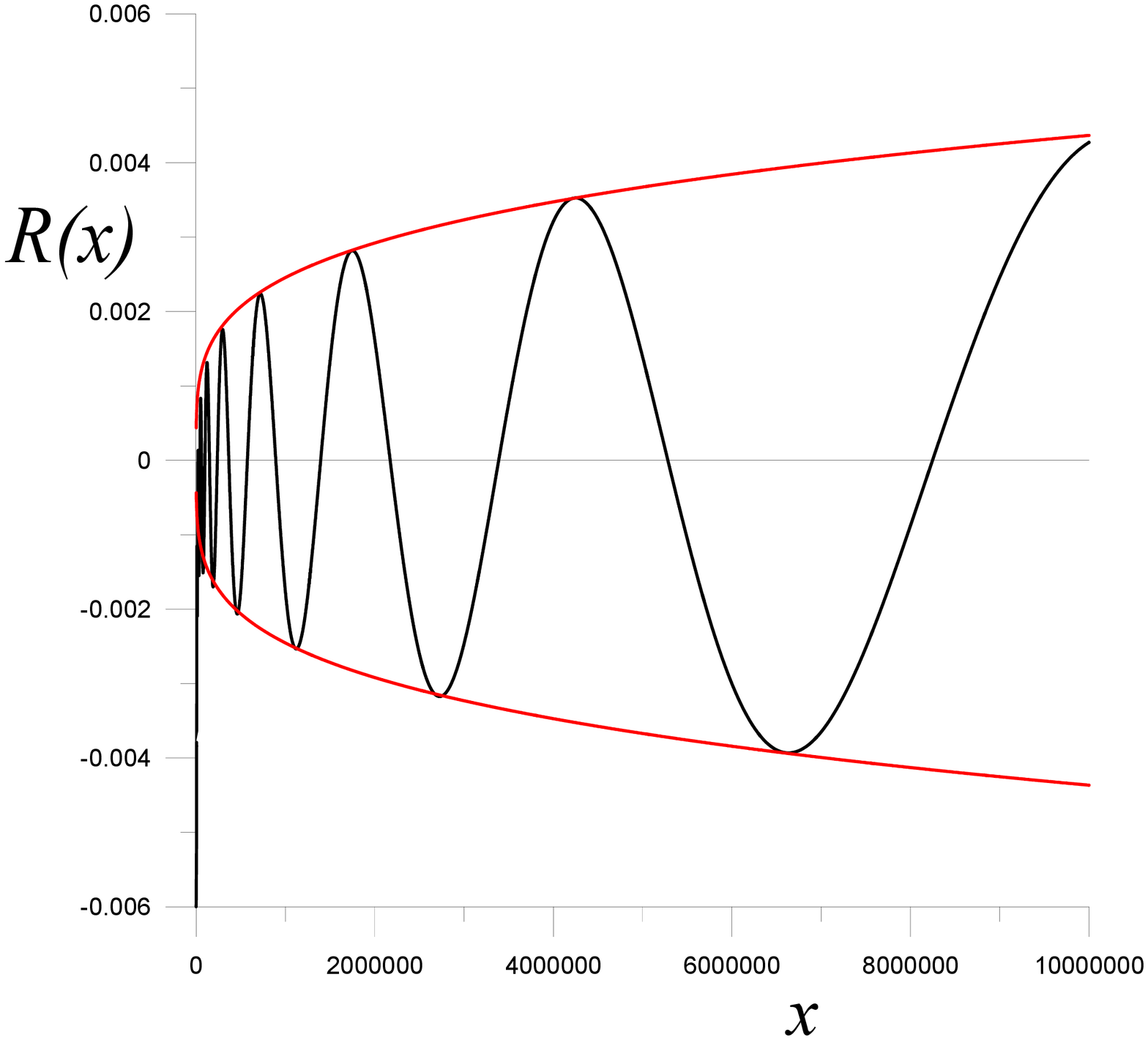} \\
\vspace{-3cm} Fig.2 The plot of $R(x)$ for $x\in(0, 10^7)$. The part of $R(x)$ smaller
than -0.006 is skipped.  \\
\end{center}
\end{minipage}
\end{figure}
\bigskip
\bigskip

\vfill

\newpage

\noindent where $\mu$ is the M\"{o}bius function:
\bee
\mu(n) \,=\,
\left\{
\begin{array}{ll}
1 & \mbox {if $ n =1 $} \\
0 & \mbox {if $n$ is divisible by a square of a prime} \\
(-1)^k & \mbox{if $n$ is a product of $k$ different primes}
\end{array}
\right.
\eee
Repeating Kummer's procedure gives:
\bee
R(x) = x \left(\frac{6}{\pi^2} + \sum_{n=1}^\infty
\frac{\mu(n)}{n^2}\left(e^{-\frac{x}{n^2}} - 1\right)\right).
\label{Riesz3}
\eee
Using this formula we were able to produce the plot of $R(x)$ for $x$ up to $10^7$, see
Fig.1 and Fig.2.  The first nontrivial zero of $R(x)$ is $x_0=1.1567116438\ldots$.
The envelops on the Fig.2 (in red) are given by the equations
\bee
y(x)=\pm A x^{\frac{1}{4}},
\eee
where $A=0.777506\ldots \times 10^{-5}$.

It is very time consuming to calculate values of the sequence $c_k$ directly
from the definition (\ref{ckmain}), see \cite{Maslanka3}, \cite{Wolf}. The point is, that for
large $j$ $\zeta(2j)$ is practically 1, and to distinguish it from 1 high precision
calculations are needed. The experience of \cite{Wolf} showed that to calculate
$c_k$ from (\ref{ckmain})  roughly $k\log_{10}(k)$ digits accuracy is needed. However in
\cite{Luis3} Baez-Duarte gave the explicit formula\footnote{There is an error
in \cite{Luis3} and there should be no minus sign in front of $c_{k-1}$ in formulae
(1.11), (1.12), (4.1), (4.11) in \cite{Luis3}.} for $c_k$ valid for large $k$:
\bee
c_{k-1}=\frac{1}{2k}\sum_{\rho} \frac{k^{\frac{\rho}{2}}\Gamma(1-\frac{\rho}{2})}{\zeta'(\rho)}
 + o(1/k)
\label{explicite}
\eee
where the sum runs for nontrivial zeros $\rho$ of $\zeta(s)$: $\zeta(\rho)=0$
and $\Im (\rho) \neq 0$.  \Mas in \cite{Maslanka3}
gives the similar formula which contains the term hidden in o(1/k) in (\ref{explicite}). 
Let us introduce the notation
\bee
\frac{\Gamma(1-\frac{\rho_i}{2})}{\zeta'(\rho_i)}= a(\rho_i)+ib(\rho_i)\equiv a_i+ ib_i.
\eee
Assuming  $\rho_i=\frac{1}{2}+i\gamma_i$ it can be shown  \cite{Wolf} that $a_i$ and  $b_i$
very  quickly decrease to zero:
\bee
\left|\frac{\Gamma(1-\frac{\rho_i}{2})}{\zeta'(\rho_i)}\right| \sim e^{-\pi\gamma_i/4}
\eee
Finally we obtain for large $k$:
\bee
c_{k-1}=\frac{1}{k^{\frac{3}{4}}}\sum_{i=1}^\infty \left\{a_i \cos\left(\frac{\gamma_i \log(k)}{2}\right) -
b_i \sin\left(\frac{\gamma_i \log(k)}{2}\right)\right\}.
\label{rownanie}
\eee
\vspace{-2.5cm}
\begin{center}
\hspace{-3truecm}
\includegraphics[width=12truecm,angle=-90]{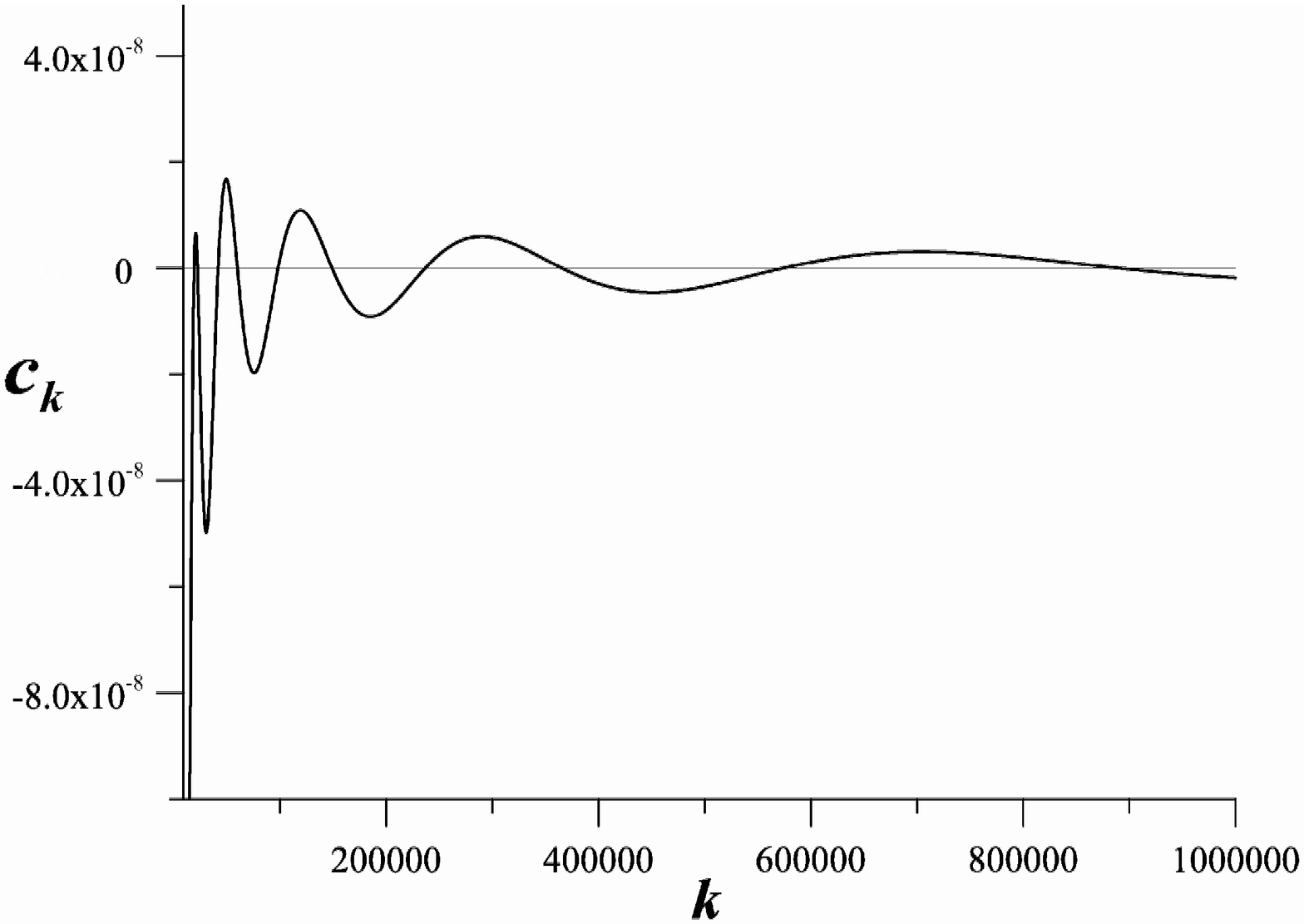}\\
\vspace{-1cm}
Fig.3 The plot of $c_k$ for $k\in(1, 10^6)$.\\
\end{center}

\bigskip

The above formula explains oscillations seen on the plots of $c_k$ published in \cite{Luis3}
and \cite{Maslanka3}, see Fig.3 . Because these curves are perfect cosine--like graphs
on the plots versus $\log(k)$ it means that in fact in the above formula (\ref{rownanie})
it suffices to maintain only  the first zero and skip all remaining  terms in
the sum.

\bigskip
\bigskip

{\bf 3. Relation between $R(x)$ and $c_k$ }

\bigskip

The values of $c_k$ can be obtained as the first elements of the sequence of forward
differences of the sequence:
\bee
f_0^0=\frac{1}{\zeta(2)} \hskip 1truecm f_1^0=\frac{1}{\zeta(4)}\hskip1truecm  f_2^0=\frac{1}{\zeta(6)}
\hskip 1truecm f_3^0=\frac{1}{\zeta(8)}\hskip 1truecm f_4^0=\frac{1}{\zeta(10)}~~~~\dots
\eee
Then we form forward differences:
\bee
f_l^k=f_l^{k-1} - f_{l+1}^{k-1}
\eee
and we have that $c_k=f_0^k$. We will recall some facts from finite difference calculus adapted for our purposes \cite{Knuth}:
Let us  define as usual  the shift operator $E$:
$$ Ef(k)=f(k+1). $$
Next we introduce sequence:
\bee
c_k=(1-E)^k f(0)=\sum_{j=0}^{k} \binom{k}{ j } (-1)^j f(j).
\label{sequence}
\eee
Then the following  equalities holds \cite{Knuth}:
$$ e^{x(1-E)}f(0)=\sum_{k=0}^{\infty}\frac{c_k}{k!} x^k, $$
$$ e^{x(1-E)}f(0)=e^xe^{-xE}f(0)=e^{x} \sum_{k=0}^{\infty}\frac{(-x)^k}{k!}f(k), $$
from which it follows that:
$$ \sum_{k=0}^{\infty}\frac{c_k}{k!} x^k= e^{x} \sum_{k=0}^{\infty}\frac{(-x)^k}{k!}f(k). $$
In our case we put
$$ f(k)=\frac{1}{\xi(2k+2)} $$
and finally we have:
\bee
\sum_{k=0}^{\infty}\frac{c_k}{k!} x^k = e^{x} \sum_{k=0}^{\infty}\frac{(-x)^k}{k!\xi(2k+2)}=
\frac{e^x}{x}R(x).
\label{rownosc1}
\eee
Thus $R(x)$ can be reconstructed from $c_k$. Vice versa, we will see later, see (\ref{bound}),
 that within some accuracy $c_k$ can be obtained from $R(x)$.
In the paper \cite{Wolf} it was suggested that the duality holds:
\bee
c_k=\mathcal{O}\left(k^{-\delta} \right) \Leftrightarrow
R(x)=\mathcal{O}\left(x^{1-\delta}\right).
\label{duality}
\eee
Putting $\delta=\frac{3}{4}-\epsilon$ gives original criteria (\ref{Riesz criterion})
and (\ref{criterion}).  In fact we will prove it in the following form:

\begin{theorem}\label{dualnosc}
The sequence $c_k$ defined by (\ref{ckmain}) decrease like $c_k=\mathcal{O}\left(k^{-\delta}
\right)$  if and only if the function $R(x)$ defined by (\ref{Riesz})
grows like $R(x)=\mathcal{O}\left(x^{1-\delta}\right)$, where $\delta<3/2$.
\end{theorem}

\medskip

{\bf Remark: } In fact $\delta$ is smaller than $3/4$, as shown by  Baez-Duarte
in \cite{Luis3}.\\

{\bf Proof:} The reasoning  that if $c_k=\mathcal{O}\left(k^{-\delta}
\right)$  then $R(x)=\mathcal{O}\left(x^{1-\delta}\right)$ we will base on the
following facts from Exercises 67 -- 71
in the Part IV of famous book of G. Polya and G. Szeg\"{o} \cite{Polya}. We summarize these
facts adapted for our purposes in the form: Let
\bee
f(x)=\sum_{k=0}^\infty a_k x^k,
\eee
where $a_k$ are positive and decrease monotonically $a_0\geq a_1 \geq a_2 \ldots
\geq a_k \geq \ldots$. Let $\alpha$ be defined by
\bee
\log a_k \sim - \frac{k\log k}{\alpha}
\eee
and next the parameter $b$ be determined from
\bee
\log f(x) \sim b x^\alpha.
\eee
Then for large $x$ the following asymptotic relation is fulfilled:
\bee
\sum_{k=1}^\infty k^\delta a_k x^k \sim  (\alpha b x^\alpha)^\delta f(x).
\eee
In our case we have $a_k = 1/k!$ thus $f(x)=e^x$ and from Stirling formula we have
$\alpha = 1$ and next $b=1$ and hence we have from above formula for large $x$:
\bee
\sum_{k=1}^\infty \frac{k^{-\delta} x^k}{k!} \sim x^{-\delta} e^x.
\eee
If we assume that $|c_k|<A k^{-\delta}$ then  we have
\bee
\sum_{k=1}^\infty \frac{|c_k| x^k}{k!} < A x^{-\delta} e^x
\eee
and  from (\ref{rownosc1}) it follows:
\bee
|R(x)| < A x^{1 - \delta}
\eee
what is a desired inequality.

We will show now the opposite implication: from $R(x)=\mathcal{O}\left(x^{1-\delta}\right)$
it follows that $c_k=\mathcal{O}\left(k^{-\delta} \right)$. In Appendix we prove
the following inequality:  
\bee
\left|\frac{R(k)}{k}-c_k\right| \leq \frac{3\sqrt{\pi}}{16}k^{-3/2} + \mathcal{O}\left(
 k^{-2}\right)
\label{bound}
\eee
Because $|c_k| -|R(k)/k|<|c_k -R(k)/k|$ and we assume $|R(k)|\leq B k^{1-\delta}$ thus
we have
\bee
|c_k|\leq Bk^{-\delta} + \mathcal{O}\left(k^{-\frac{3}{2}} \right)
\eee
To avoid nonsense $\delta$ should be smaller than 3/2 and in fact Baez-Duarte showed
\cite{Luis3} that existence of zeros on the critical axis requires $\delta<3/4$.

\hfill $\square$\\

The comparison of the above
bound (\ref{bound}) with real computer data is given in the Fig.4. Here the fit (red line)
was obtained
by the least square method  from the data  with $k>10000$ to avoid transient regime  and it is
given by the equation $y=0.0117483x^{-1.52655}$. The fact that approximately $c_k \approx R(k)/k$
was observed previously by S. Beltraminelli and  D. Merlini \cite{Merlini}.
It can be explained heuristically as follows: Baez-Duarte gives in \cite{Luis3} despite
(\ref{ckmain}) a few  formulae for $c_k$. We need
here the following expression being the transformation of (\ref{ckmain}):
\bee
c_k=\sum_{n=1}^{\infty} \frac{\mu(n)}{n^2}\left(1-\frac{1}{n^2}\right)^k.
\label{c k moebius}
\eee
For large $k$ we can write:
\bee
c_k=\sum_{n=1}^{\infty} \frac{\mu(n)}{n^2}\left(1-\frac{k}{kn^2}\right)^k \approx
\sum_{n=1}^{\infty} \frac{\mu(n)}{n^2}e^{-n^2/k}
\eee
and comparing it with (\ref{Riesz2}) we get $c_k\approx R(k)/k$ for large $k$.\\

\bigskip

\begin{minipage}{12.8cm}
\begin{center}
\vspace{-2.0cm}
\includegraphics[width=11cm ]{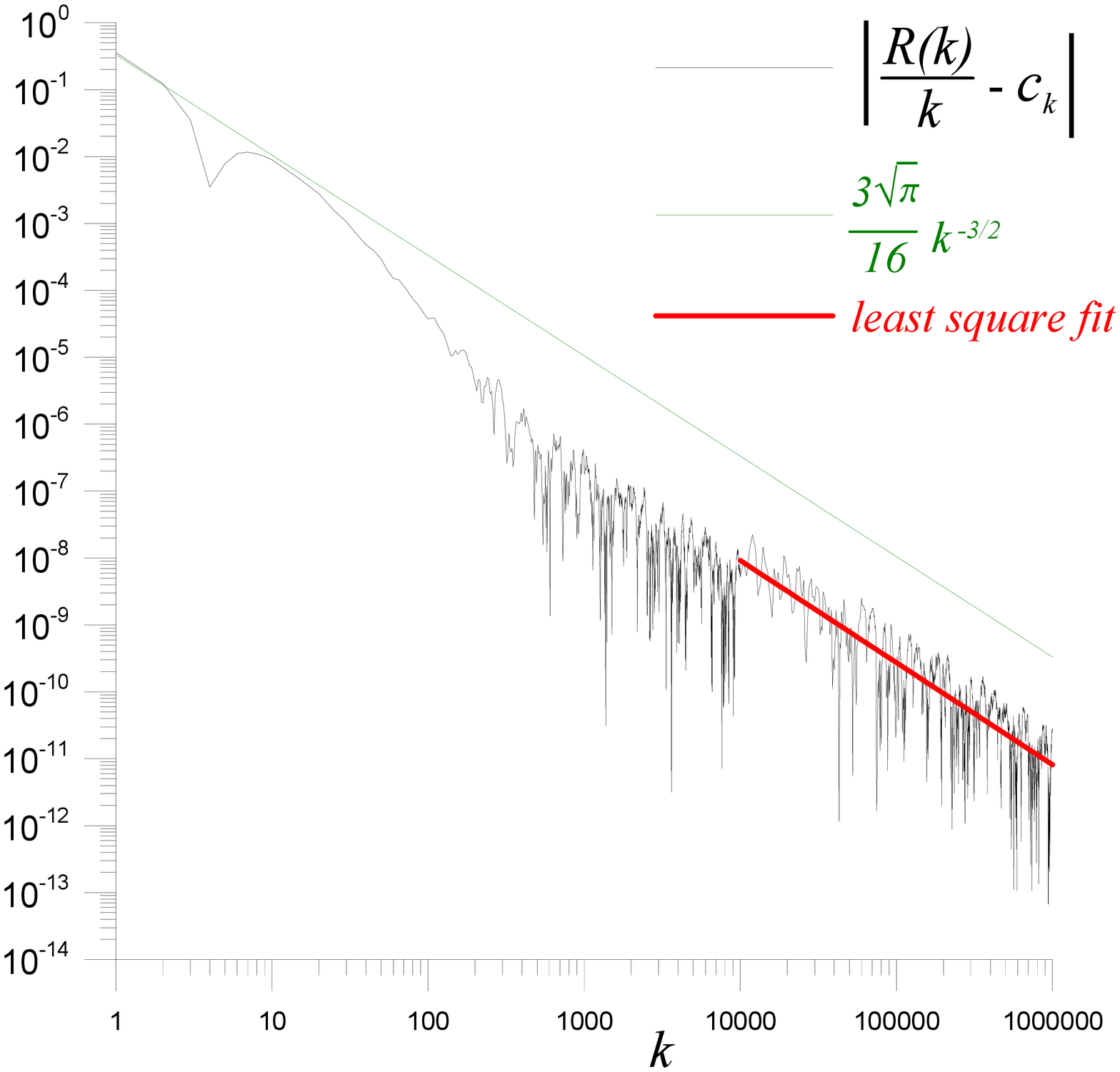}\\
\vspace{-2.5cm}
Fig. 4 The log-log plot of $|c_k -R(k)/k|$ for $k\in(0,10^6)$.
\end{center}
\end{minipage}\\


{\bf 4 Some other relations}\\

Using the formula (\ref{c k moebius})
it is possible to calculate the alternating sum of $c_k$:
\bee
 \sum_{k=0}^{\infty}(-1)^k c_k  = \sum_{k=1}^{\infty}\frac{1}{2^k} \frac{1}{\zeta(2k)}.
\label{suma}
\eee
Numerically this sum is $\sum_{k=0}^{\infty}(-1)^k c_k = 0.782527985325384234576688\ldots$.
This number probably can not be expressed by other known constants, because the Simon Plouffe
inverter failed to find any relation \cite{Plouffe}.    By the Abel's summation the r.h.s.
can be written as:
\bee
\sum_{k=1}^{\infty}\frac{1}{2^k} \frac{1}{\zeta(2k)} =
1 + \int_2^\infty \left(1-\frac{1}{2^{\lfloor x/2 \rfloor }} \right) \frac{\zeta'(x)}{\zeta^2(x)} dx.
\eee
In fact more general than (\ref{suma}) formula holds:
\bee
\sum_{k=0}^{\infty}c_k  s^k = \frac{1}{1-s}
\sum_{k=0}^{\infty}\left( \frac{-s}{1-s}\right)^k \frac{1}{\zeta(2k+2)},
\label{suma_2}
\eee
where $-1 \leq s <\frac{1}{2}$. Here we have made use of the identity
\bee
\sum_{n=1}^{\infty} \left(1-\frac{1}{n^2}\right)^k s^k =
\frac{1}{1-s}\sum_{k=0}^{\infty}\left( \frac{-s}{1-s}\right)^k \frac{1}{n^{2k}}.
\eee
The l.h.s. is convergent for $-1\leq s <1$ while the r.h.s. converges for  $-\infty<  s<1/2$.

The question of the convergence of the sum $\sum_{k=0}^\infty c_k $
is much more complicated. Formally summing both sides of (\ref{sequence}) we get:
\bee
\sum_{k=0}^\infty c_k = E^{-1}f(0) = f(-1)
\label{formal}
\eee
As in our case $f(k)=1/\zeta(2k+2)$ we have
\bee
\sum_{k=0}^\infty c_k = \frac{1}{\zeta(0)} = -2
\eee
because $\zeta(0)=-\frac{1}{2}$, see e.g. \cite{Titchmarsh}, p.19. The partial sums
$\sum_{k=0}^n c_k $ indeed initially tend from above to -2, but for $n\approx 91000$
the partial sum crosses -2 and around $n\approx 100000$ the partial sum starts to increase.
These oscillations begins to repeat with growing amplitude around -2. The Fig. 5
shows  the plot of  distances of the partial sums $\sum_{k=0}^n c_k $ from -2. Let us remark
that at $n\sim 10^8$ the amplitude is rather very small: of the order 0.001.
When we retain in (\ref{rownanie}) only the first zero $\gamma_1$
it can be shown that this amplitude grows like $n^{1/4}$, thus it appears that the above
formal derivation (\ref{formal}) is wrong and the sum $\sum_{k=0}^n  c_k$ is divergent.

\begin{minipage}{12.8cm}
\begin{center}
\vspace{-2.0cm}
\includegraphics[width=10cm, scale=0.8]{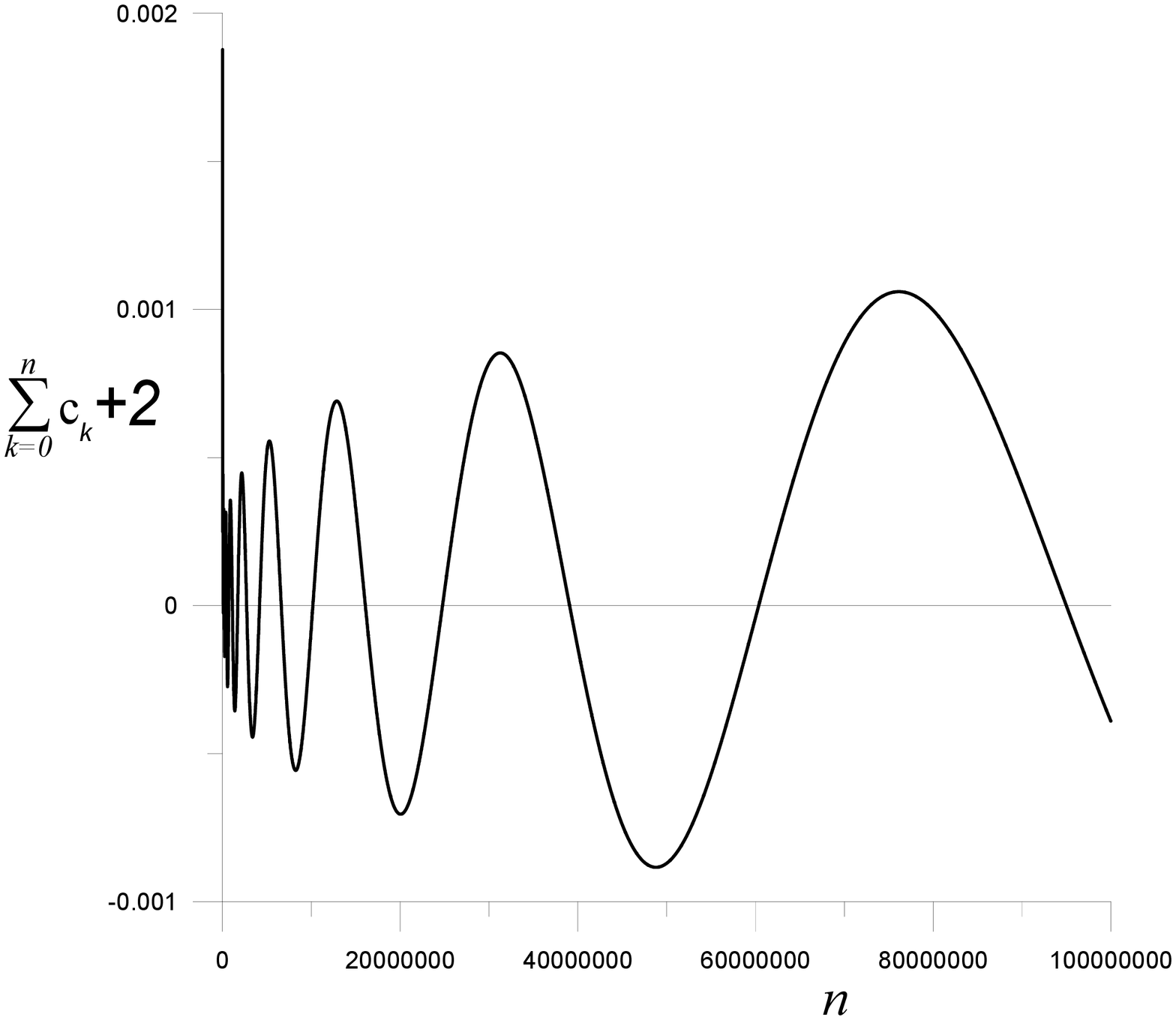}\\
\vspace{-2.5cm}
Fig. 5 The distance from -2 of the partial sums $\sum_{k=0}^n$ for $n=1,\ldots 10^8$.
\end{center}
\end{minipage}

\bigskip

We have made analogous plot of
partial sums $\sum_{k=0}^n (-1)^k c_k$ and there we have seen oscillations around
the limit value $s=0.782527985 \ldots$  of {\it decreasing} amplitude. Thus we speculate,
that this partial sums behave as  $c_k$ and $R(x)$ accordingly:
\bee
\left| \sum_{k=0}^n (-1)^k c_k - s\right | = \mathcal{O}\left(n^{-\frac{3}{4}}\right),
\eee
\bee
\left| \sum_{k=0}^n  c_k +2 \right | = \mathcal{O}\left(n^{\frac{1}{4}}\right).
\eee

\bigskip

Finally we would like to argue in favour of the two strange approximate equalities.
Both follows from  $c_k \approx R(k)/k$ for large $k$. The first follows when write this
relation with the help of (\ref{Riesz}) and (\ref{ckmain}):
\bee
\sum_{j=0}^\infty \frac{(-1)^{j}k^{j}}{j!\zeta(2j+2)}\approx
\sum_{j=0}^k (-1)^j \binom{k}{j}\frac{1}{\zeta(2j+2)}.
\eee
On both sides there appears inverses of $\zeta(2n)$. We have checked numerically that the
difference between these two sums very quickly tends to zero.

The second formula we get when in (\ref{rownosc1}) we put instead of $c_k$ simply  $R(k)/k$:
\bee
 \frac{e^x}{x}R(x) \approx \sum_{k=0}^\infty R(k) \frac{x^k}{k \ k!}
\eee
Thus we get  $R(x)$ as ``entangled'' combination of $R(k)$ at positive integers.
We end asking the question:  Will such
a kind of constraint help to prove (\ref{Riesz criterion})?

\bigskip

{\bf Acknowledgement} We thank Prof. L. Baez-Duarte and Prof. K. \Mas for e-mail exchange.
To prepare  data for some figures we have used  the free package PARI/GP \cite{Pari}.

\bigskip

{\bf Appendix }

In this appendix we will calculate the error of the  approximation $c_k\approx R(k)/k$.
Looking at (\ref{Riesz2}) and (\ref{c k moebius}) we see that we have to estimate the sum:
\bee
\left|\sum_{n=1}^\infty \frac{\mu(n)}{n^2} ~e^{-\frac{k}{n^2}} - \sum_{n=1}^{\infty} \frac{\mu(n)}{n^2}\left(1-\frac{1}{n^2}\right)^k\right|<
\sum_{n=1}^\infty  \left| \frac{1}{n^2} ~e^{-\frac{k}{n^2}} - \frac{1}{n^2}\left(1-\frac{1}{n^2}\right)^k\right|
\label{Start}
\eee
Instead of $\mu(n)$ we have put 1. Let $h(x)$ denote for $1\leq x$:
\be
h(x)=\frac{1}{x^2}\exp(-k/x^2)-\frac{1}{x^2}\left(1-\frac{1}{x^2}\right)^k.
\ee
This function is bounded by:
\be
0 < h(x) \leq \frac{27}{2e^3k^2}+\frac{128}{e^4k^3} .
\ee
and has one maximum. Thus we can apply the rule:
\be
\sum_{n=1}^{\infty} h(n) \leq \max_n h(x) + \int_{1}^{\infty} h(x)dx.
\ee
The integral is estimated as
$$ \int_1^{\infty} \frac{1}{x^2}\left(\exp(-k/x^2)-(1-1/x^2)^k\right)dx=\int_0^1 (e^{-ky^2}-(1-y^2)^k)dy. $$
$$ \int_0^1 e^{-ky^2}dy=k^{-1/2}\int_0^{\sqrt{k}} e^{-y^2}dy \leq k^{-1/2} \int_0^{\infty} e^{-y^2}dy
= \frac{\sqrt{\pi}}{2\sqrt{k}}. $$
\bee
 \int_0^1 (1-y^2)^k dy = \frac{4^k}{\binom{2k}{k}(2k+1)}
\geq \frac{ \sqrt{\pi k}}{2k+1}\left(1+\frac{1}{8k}-\frac{1}{72k^2}\right).
\label{Beta}
\eee
Here the Stirling formula in the form
\be
k!=\sqrt{2\pi k} \; k^k e^{-k+\theta(k)}, \quad \frac{1}{12k+1} < \theta(k) < \frac{1}{12k}
\ee
was used. Collecting all above estimations we obtain:
$$ \int_1^{\infty} \frac{1}{x^2}(\exp(-k/x^2)-(1-1/x^2)^k)dx \leq $$
$$ \leq \frac{\sqrt{\pi}}{2\sqrt{k}}- \frac{ \sqrt{\pi k}}{2k+1}\left(1+\frac{1}{8k}-\frac{1}{72k^2}\right)
<\frac{3\sqrt{\pi}}{16}k^{-3/2}+\frac{\sqrt{\pi}}{144}k^{-5/2}. $$
and finally from the starting sum (\ref{Start}) we get the desired inequality:
$$ | R(k)/k - c_k |\leq \frac{3\sqrt{\pi}}{16}k^{-3/2} + \frac{27}{2}e^{-3}k^{-2}+\frac{\sqrt{\pi}}{144}k^{-5/2}
+128e^{-4}k^{-3}. $$
For $k>16$ it suffices to retain  in the above inequality on the r.h.s only the leading
term $k^{-3/2}$.
Let us remark that the integral (\ref{Beta}) can be also taken from tables
as it is the Euler Beta integral:
\be
\int_0^1(1-y^2)^k dy = \frac{1}{2}B\left(\frac{1}{2}, k+1\right)
\ee
and from
\be
B(a,x)\sim x^{-a} \Gamma(a) ~~~~~~{\rm for}~~x~~{\rm large}
\ee
see e.g. \cite{Titchmarsh2}\S 1.8.7,  we get:
\be
\int_0^1(1-y^2)^k dy \sim \frac{\Gamma(\frac{1}{2})}{2\sqrt{k+1}} =\frac{\sqrt{\pi}}{2\sqrt{k+1}}
\ee
what for large $k$ reproduces leading term in (\ref{Beta}).

\end{document}